\documentclass[10pt,reqno]{article}
\usepackage{geometry}
\usepackage{amsmath}
\usepackage{amsthm}
\usepackage{amssymb}
\usepackage{prodint}
\usepackage{color}
\usepackage{authblk}
\usepackage[round]{natbib}
\usepackage{rotating}
\usepackage{subfigure}
\usepackage{subfig}
\usepackage{times}
\usepackage{epsfig}
\usepackage{amsfonts}
\usepackage{float}
\usepackage{bbm}
\usepackage{latexsym}
\usepackage{graphicx}
\usepackage{caption}
\usepackage{booktabs}
\usepackage{mathabx}
\usepackage{xcolor}

\newtheorem{thm}{Theorem}
\newtheorem{lemma}{Lemma}
\newtheorem{rem}{Remark}
\newtheorem{cor}{Corollary}

\newcommand{\rnc}{\renewcommand}
\newcommand{\nc}{\newcommand}
\newcommand{\mrm}{\mathrm}
\renewcommand{\b}{\textbf}

\newcommand{\bs}{\boldsymbol}
\nc{\mb}{\mathbb}
\nc{\mc}{\mathcal}
\nc{\E}{\mb{E}}
\nc{\N}{\mb{N}}
\nc{\R}{\mb{R}}
\nc{\Q}{\mb{Q}}
\rnc{\P}{\mrm P}
\nc{\bP}{\b P}
\nc{\bA}{\b A}
\rnc{\d}{\mrm d}
\nc{\C}{\mc{C}}
\nc{\D}{\mc{D}}
\nc{\B}{\mc{B}}
\nc{\gDg}{\stackrel{d}{=}}
\nc{\oPo}{\stackrel{\mrm p}{\longrightarrow}}
\nc{\oWo}{\stackrel{w}{\longrightarrow}}
\nc{\oDo}{\stackrel{d}{\longrightarrow}}
\nc{\nae}{Nelson-Aalen estimator}
\nc{\aje}{Aalen-Johansen estimator}
\nc{\naeL}{Nelson-Aalen estimator\ }
\nc{\ajeL}{Aalen-Johansen estimator\ }
\nc{\CIF}{cumulative incidence function}
\nc{\CIFL}{cumulative incidence function\ }
\nc{\wh}{\widehat}
\nc{\leb}{\lambda \!\! \lambda}
\nc{\red}{\color{red}}
\nc{\tcp}{\textcolor{purple}}

\newcommand\blfootnote[1]{%
  \begingroup
  \renewcommand\thefootnote{}\footnote{#1}%
  \addtocounter{footnote}{-1}%
  \endgroup
} 

\begin{document}

\title{\Large \bf Survival of the Fittest Group: \\ Factorial Analyses of Treatment Effects \\ under Independent Right-Censoring}
\author[1$*$]{Dennis Dobler}
\author[2]{Markus Pauly}

\affil[1]{Department of Mathematics, Vrije Universiteit Amsterdam, Netherlands.}
\affil[2]{Institute of Statistics, Ulm University, Germany.}

\maketitle

\begin{abstract}
\noindent 
\blfootnote{${}^*$ e-mail:  d.dobler@vu.nl}

\noindent This paper introduces new effect parameters for factorial survival designs with possibly right-censored time-to-event data. In the special case of a two-sample design it coincides with the concordance or Wilcoxon parameter in survival analysis. More generally, the new parameters describe treatment or interaction effects and we develop estimates and tests to infer their presence. We rigorously study the asymptotic properties by means of empirical process techniques and additionally suggest wild bootstrapping for a consistent and distribution-free application of the inference procedures. 
The small sample performance is discussed based on simulation results.
The practical usefulness of the developed methodology is exemplified on a data example about patients with colon cancer by conducting one- and two-factorial analyses.
\end{abstract}


\vfill
\vfill


\newpage

\section{Motivation and Introduction}
Many biomedical and clinical trials are planned as factorial designs. Here, not only the (main) effects of separate factors but also interaction effects that are related to possibly complex factor combinations are of importance. 
Such interaction effects
may even alter the interpretation of main effects leading to the established comment by \cite{lubsen1994factorial} that 
{\it`it is desirable for reports of factorial trials to include estimates of the interaction between the treatments'}.
 
On the other hand, {\it nonparametric} estimation and the inference of adequate effects in such designs can be rather involved. 
In particular, most existing inference procedures have focused on testing hypotheses formulated in terms of distribution functions \citep{brunner97, brunner2001nonparametric, gao2005unified, gao2008nonparametric, gao2008nonparametriconeway, akritas2011nonparametric, dutta16, friedrich2016CSDA, UKP2016} which cannot be inverted to obtain confidence intervals or regions for meaningful effects. 
Only recently, nonparametric methods for inferring adequate effects in general factorial designs with independent and dependent observations have been established \citep{rankFD, brunnerrank2016, umlauft2017wild, dobler2017nonparametric2}. These procedures are, however, only developed for completely observed data and not applicable for 
partially observed time-to-event data. 
Since many clinical studies are concerned with survival outcomes, adequate statistical inference methods for complex factorial time-to-event designs are of particular interest. 

To detect main effects, weighted logrank tests or their extensions may be applied in case of two or multiple samples  
\citep{mantel1966evaluation, abgk93, ehm1995power, liu1995design, janssen1997two, bathke2009combined, yang2010improved, fleming2011counting, Brendel_etal_2014}. However, these procedures only infer conclusions
in terms of cumulative hazard functions and cannot be applied to obtain 
concrete {\it effect parameters} with informative confidence intervals nor tests for the presence of interactions. 
In practice, interaction effects are usually modeled with the help of
Cox-, Aalen- or even Cox-Aalen regression models \citep{cox72, scheike2002additive, scheike2003extensions} with factors as covariates and incorporated interaction terms. However, although very flexible, these models are usually more driven towards hazards modeling by continuous covariates while the incorporation of several factor variables (e.g., via multiple dummy variables per factor) can become cumbersome; especially when interactions are incorporated, see also \citet{green2002factorial} and \cite{crowley2012} for the uncensored case. 

The above problems directly motivate a nonparametric approach for estimating and inferring main and interaction effects in factorial designs with censored observations. 
So far, the only existing methods in this context are given by the nonparametric survival procedures of \cite{akritas97} and \cite{akritas2011nonparametric}.
They are based on a purely nonparametric model that does not require any multiplicative or additive structure of the hazards and can even be applied for arbitrary, possibly non-continuous survival distributions (i.e., it can be readily used for survival times rounded to days, weeks or months). 
Moreover, it leads to tests for
main and interaction effects in case of independent right-censored data. 
However, these tests suffer from several drawbacks: 
the procedure is based on a rather strong assumption on the underlying censorship distribution 
which is often hard to verify in practical situations. In addition, 
null hypotheses are only formulated in terms of distribution functions. 
As a result, there is no direct quantification and estimation of main and interaction effects in terms of confidence intervals as, e.g., required by regulatory authorities (ICH E9 Guideline, 1998, p. 25).

This is to be changed in the current paper. We develop and rigorously analyze nonparametric inference procedures, i.e. tests and confidence intervals, for meaningful effect sizes in factorial survival designs, where data may be subject to random right-censoring. 

Similar to the adaption of the \citet{brunner2000nonparametric} test to the two-sample survival set-up by \citet{dobler2016bootstrap}, we consider the recently proposed unweighted nonparametric effects of \citet{brunnerrank2016} and extend their ansatz to a general survival setting. In the special case of proportional hazards, these effects have a direct relationship to hazard ratios in two-sample settings \citep{bruckner2017sequential} while they remain meaningful in case of non-proportional hazards. This fact makes the effect sizes even more appealing for practical purposes.

The paper is organized as follows. The statistical model and important results on the basic estimators are presented in Section~\ref{sec:mod}. The resulting test statistic for the null hypotheses of interest is stated and mathematically analyzed in Section~\ref{sec:test_stat}. Since the asymptotic distribution of the test statistic depends on unknown parameters, 
we propose a distribution-free multiplier resampling approach in Section~\ref{sec:wbs} and prove its consistency. 
In Section~\ref{sec:simus} it is supplemented by a simulation study to assess the finite sample properties of the proposed procedure.
They are then exemplified on a colon cancer study in Section~\ref{sec:data_Example}, where in the original study~\citep{moertel90} the analysis was made in terms of Cox models.
Finally, the paper closes with concluding comments in Section~\ref{sec:dis}.
All proofs are deferred to the technical Appendix.

\section{The set-up}
\label{sec:mod}
To establish the general model, we consider sequences of mutually independent random variables
\begin{equation}\label{eq:mod}
 T_{ik} \stackrel{\text{ind}}{\sim} S_i \quad{\mbox{and}}\quad C_{ik} \stackrel{\text{ind}}{\sim} G_i \qquad( i = 1, \dots, d, \ k = 1, \dots, n_i),
\end{equation}
where $T_{ik}$ denotes the actual survival time of subject $k$ in group $i$ and $C_{ik}$ the corresponding censoring variable. Moreover, to even allow for ties or survival times rounded to weeks or months, the survival functions $S_i$ and $G_i$, $i=1, \dots, d$, defined on $(0,\infty)$ may be {\it possibly discontinuous}.
That is, the corresponding hazard rates may, but need not exist. 
The actually observable data consist of the right-censored survival times $X_{ik} = T_{ik} \wedge C_{ik}$ and the uncensoring indicators $\delta_{ik} = 1 \{ T_{ik} \leq C_{ik} \}$, 
 $i = 1, \dots, d, \ k = 1, \dots, n_i$. In this set-up, a factorial structure can be incorporated by splitting up indices, see Section~\ref{sec:simus} for details. 
 
 In the special case of $d=2$ groups with continuous survival times \citet{efron67} introduced an estimator for the {\it concordance probability}
 $$
  w= P(T_{11} > T_{21}) = - \int  S_1 dS_2
 $$
 that a randomly chosen subject from the first group survives longer than someone from the second group. 
 If all subjects are completely observable, this effect size $w$ reduces to the well-known Mann-Whitney effect underlying the \citet{brunner2000nonparametric} test.
 Inference procedures for $w$ and related quantities in survival set-ups (such as the concordance parameter or the average hazard ratio)  have, e.g. been developed by \citet{bruckner2017sequential, dobler2016bootstrap}. However, an extension of the definition of $w$ to the more general design \eqref{eq:mod}, allowing for an arbitrary factorial structure, is not straightforward. In particular, for the case of completely observed data, \citet{brunnerrank2016} and \citet{brunner2018ranks} point out several pitfalls that may lead to paradoxical results when working with a `wrong' extension of $w$. Adopting their solution to the present situation, we  introduce an additional `benchmark' survival time $Z$, independent of the above, with averaged survival function $Z\sim \bar S = \frac1d \sum_{i=1}^d S_i$. This is used to extend $w$ to 
 $$\tilde p_i = P(T_{i1} > Z) + \frac12 P(T_{i1} = Z) = 
 - \int S_i^\pm \d \bar S, $$
 where the superscript $\pm$ denotes the average of a right-continuous function and its left-continuous version. 
 The use of such normalized survival functions adequately handles discrete components of the survival distribution, i.e.\ ties in the data are explicitly allowed.
  
 The choice of the effect parameter $\tilde p_i$ is motivated by recent findings on nonparametric analyses of factorial designs with 
 complete observations in \cite{brunnerrank2016, brunner2018ranks}. \nocite{dobler2017nonparametric2} They stress that other choices, e.g.
 pairwise-comparisons of all concordance probabilities $w$ or comparisons with the weighted survival function
$\sum_{i=1}^d \frac{n_i}{N}S_i$ instead of $\bar S$, may easily result in paradoxical outcomes. This is no issue for the effects $\tilde p_i$ which are sample size independent. 
 For later calculations, we emphasize that the effect parameters are balanced in the mean.
 In particular, we have 
 $$
 \frac1d \sum_{i=1}^d \tilde p_i = - \int \bar S^\pm \d \bar S = \frac12
\quad\text{and}\quad
 \tilde p_i = - \sum_{\substack{j = 1 \\ j \neq i}}^d \int S_i^\pm \d S_j + \frac1{2d}.
 $$
 From a practical point of view, estimation of the $\tilde p_i$'s would need 'arbitrarily' large survival times since the integral is defined on $(0,\infty)$.
 However,  every study ends at a certain point in time. For practical applicability, we therefore assume that the censoring times are bounded and we have to modify the $\tilde p_i$'s accordingly:  
 denote by $\tau > 0$ the largest possible censoring time. 
 In comparisons of survival times, which belong to different groups and which exceed $\tau$,
 no group shall be favored.  In other words, the remaining mass has to be split up equally among the groups.  Technically, this is realized by setting the remaining mass of the survival functions to zero: $S_i(\tau) =0$. 
 Redefining $S_i$ and $\bar S$ from now on as the survival functions of $\min(T_{i1}, \tau)$ and $\min(Z, \tau)$, respectively, this translates into the {\it nonparametric concordance effects} 
 \begin{equation}
 p_i = P(\min(T_{i1}, \tau) > \min(Z, \tau)) + \frac12 P(\min(T_{i1}, \tau) = \min(Z, \tau)) = - \int S_i^\pm \d \bar S.
 \end{equation}
 Obviously, all of the above-discussed positive properties of the effects parameter $\tilde p_i$ also transfer to the nonparametric concordance effects $p_i$:
 it is a meaningful effect measure for ordinal and metric data, sample size independent, and allows for a suitable treatment of ties. 
 
  We aggregate all effects into the vector $\b p = (p_1, \dots, p_d)'$ and borrow a trick from \cite{konietschke2012rank} and \cite{brunnerrank2016} to express them as
  \begin{equation}\label{eq:p_as_w_expression}
\b p  = \Big(\b I_d \otimes \frac1d \b 1_d'\Big) \cdot (\b w_1', \dots,  \b w_d')' =: \b E_d \cdot \b w.   
  \end{equation}
Here, $\b w_i = (w_{1i}, \dots, w_{di})' = - \int S_i^\pm \d \b S$ is the $\R^d$-vector of effects for direct comparisons of group $i$ with respect to all groups $j=1,\dots,d$,
 and $\b S = (S_1, \dots, S_d)'$ is the aggregation of all survival functions. Moreover, $\b I_d$ denotes the identity matrix in $\R^d$, $\b 1_d$  the $d$-dimensional vector of 1's and the symbol $\otimes$ denotes the Kronecker product. In this way the $i$th entry of $\b w_i$ is $w_{ii} = \frac12$ which makes sense because equal groups should be valued equally high. 
 Anyhow, Equation~\eqref{eq:p_as_w_expression} shows that the problem of estimating $\b p$ reduces to the estimation of the pair-wise effects $w_{ji}$. But this can be achieved by substituting each involved survival function $S_i$ by its Kaplan-Meier estimator $\wh S_i, i=1,\dots,d$ \citep{kaplan58}. Proceeding in this way we denote by 
 $\wh {\b w}$ and $\wh {\b w}_i$ these estimated counterparts of $\b w$ and $\b w_i$. Let $N = \sum_{i=1}^d n_i$ be the total sample size. 
 Below we establish the asymptotic normality of 
 $\sqrt{N}(\wh {\b w} - {\b w})$ under the following framework:
 \begin{align}
 \label{eq:sample_size_conv}
    N^{-1} \bs n :=  \Big(\frac{n_1}{N}, \dots, \frac{n_d}{N}\Big)' \rightarrow \bs \lambda := (\lambda_1, \dots, \lambda_d)' \in (0,1)^d
 \end{align}
 as $\min \bs n \rightarrow \infty$. To give a detailed description of the resulting asymptotic covariance structure, however, we first have to introduce some additional notation: 
  Let $D[0,\tau]$ be the space of all c\`adl\`ag-functions on $[0,\tau]$, equipped with the Skorokhod metric, and $BV[0,\tau] \subset D[0,\tau]$ its subspace of c\`adl\`ag-functions with bounded variation.
 For the subsequent arguments it is essential that we can represent  $\b w = \phi \circ \b S $  as a functional of $\b S$. In particular, the functional
 $$ \phi: \  (BV[0,\tau])^d \rightarrow \R^{d^2}, \qquad (f_1, \dots, f_d)' \longmapsto \Big( - \int f_i^\pm \d f_j \Big)_{i,j = 1}^d, $$
 with inner index $j$, is Hadamard-differentiable at $\b S$; see the proof of Lemma~\ref{lem:w} below for details.
 We denote its Hadamard-derivative at $\b S$ by $\d \phi_{\b S}$, which is a continuous linear functional.
 For technical reasons, we assume throughout that $P(T_{i1} > \tau) >0$ for all groups $i=1, \dots, d$.
 We may now state the first preliminary but essential convergence result.
 \begin{lemma}\label{lem:w}
 Under the asymptotic regime \eqref{eq:sample_size_conv} we have
 $$  \sqrt N ( \wh {\b w} - \b w ) \oDo \b W, $$
 where $\b W$ has a centered multivariate normal distribution on $\R^{d^2}$.
 \end{lemma}
 In particular, we can write $\b W = \d \phi_{\b S} \cdot diag(\bs \lambda)^{-1/2} \b U$,
 where $\b U$ consists of independent, zero-mean Gaussian processes $U_1, \dots, U_d$ with covariance functions 
 $$\Gamma_i(r,s) = S_i(r) S_i(s) \int_0^{r \wedge s} \frac{\d \Lambda_i}{S_{i-} G_{i-} (1 - \Delta \Lambda_i)}, \quad i=1,\dots,d,$$
 where $\Lambda_i$ denotes the cumulative hazard function corresponding to $S_i = \prodi (1 - \d \Lambda_i)$, $i=1, \dots, d$; 
 the symbol $\prodi$ denotes the product integral \citep{gill90}. 
 Here, a minus sign in a subscript indicates the left-continuous version of a function and $\Delta \Lambda = \Lambda - \Lambda_-$ is the jump size function of $\Lambda$.
Note that the covariance matrix of $\b W$ is singular; in particular, $ ( \d \phi_{\b S} \cdot  diag(\bs \lambda)^{-1/2} \b U )_{i,i} = 0 $ for all $i = 1, \dots, d$.
 The other entries ($i \neq j$) are distributed as follows:
 \begin{align*}
  ( \d \phi_{\b S} \cdot  diag(\bs \lambda)^{-1/2}  \b U )_{i,j} 
  =  \frac{1}{\sqrt{\lambda_i}} \int U_i^\pm \d S_j 
    - \frac{1}{\sqrt{\lambda_j}} \int U_j^\pm \d S_i  
  \sim N \Big(0, \frac{1}{\lambda_i} \int \int \Gamma_i^{\pm \pm} \d S_j \d S_j 
  + \frac{1}{\lambda_j} \int \int \Gamma_j^{\pm \pm} \d S_i \d S_i \Big) .
 \end{align*}
 Here, the double appearance of $\pm$ signs means the average of all four combinations of left- and right- continuous versions in both arguments of a two-parameter function.

 Let us now turn to the estimation of the nonparametric concordance effects $\b p$.
 A matrix multiplication of $\wh w$ with $\b E_d$ from the left is basically the same as taking the mean with respect to the inner index $j$. This immediately brings us to the first main result:
 \begin{thm}
  \label{thm:p}
  Under the asymptotic regime \eqref{eq:sample_size_conv} we have
 $$\sqrt N ( \wh {\b p} - \b p) := \sqrt N \b E_d ( \wh {\b w} - \b w ) \oDo \b E_d \b W =\Big( \frac1d \sum_{i=1}^d \frac{1}{\sqrt{\lambda_i}} \int U_i^\pm \d S_j - \frac{1}{\sqrt{\lambda_j}} \int U_j^\pm \d \bar S  \Big)_{j=1}^d,$$ 
 where $ \b E_d \b W$ has the variance-covariance matrix $\b V$ with the following entries:
 $$ V_{ii} = \frac{1}{\lambda_i} \int \int \Gamma_i^{\pm \pm} \d \bar S \d \Big( \bar S - \frac2d S_i \Big)
  + \frac{1}{d^2} \sum_{j=1}^d \frac{1}{\lambda_j} \int \int \Gamma_j^{\pm \pm} \d S_i \d S_i $$
  in the $i$th diagonal entry, $i = 1, \dots, d$, and 
  $$ V_{ij} = \frac{1}{d^2} \sum_{j=1}^d \frac1{\lambda_j} \int \int \Gamma_j^{\pm \pm} \d S_i \d S_j
   - \frac1d \frac1{\lambda_i} \int \int \Gamma_i^{\pm \pm}  \d \bar S \d S_j 
   - \frac1d \frac1{\lambda_j} \int \int \Gamma_j^{\pm \pm}  \d \bar S \d S_i $$
   in the off-diagonal entries $(i,j)$, \ $i \neq j$.
 \end{thm}
 A more compact form of the matrix $\b V$ is given in Appendix~\ref{app:asy.cov}.

  \section{Choice of Test Statistic}
  \label{sec:test_stat}
  
  In order to develop hypothesis tests based on the estimator $\wh{\b p}$, 
  we next need to find a consistent estimator $\wh {\b V}_N$ for $\b V$. 
  A natural choice  is to plug in estimators for all unknown quantities that are involved in $\b V$.
  In particular, we use the Kaplan-Meier estimators for all survival functions and
  $\wh \Gamma_i(s,t) = \wh S_i(s) \wh S_i(t) n_i \int_0^{s \wedge t} [ Y_i (1 - \Delta \wh \Lambda_i)]^{-1} \d \wh \Lambda_i$ for each covariance function $\Gamma_i$,
  where $Y_i$ is the number at risk process and $\wh \Lambda_i$ is the Nelson-Aalen estimator for the cumulative hazard matrix in group $i$.
  Note that if $\Delta \wh \Lambda_i(u) = 1 $, we also have $\wh S_i(u) = 0$ in which case we let $\wh \Gamma_i(s,t) = 0$ if $s \geq u$ or $t \geq u$.
  We denote the resulting covariance matrix estimator by $\wh {\b V}_N$.

  \begin{lemma}
   \label{lem:gamma}
  Under the asymptotic regime \eqref{eq:sample_size_conv} we have the consistency
  $\wh {\b V}_N \oPo \b V$. 
  \end{lemma}

  All of the developed convergence results are now utilized to find the most natural test statistic.
  First note that the asymptotic covariance matrix $\b V$ is singular 
  since $\b 1_d' \sqrt{N} (\wh {\b p} - \b p) \equiv 0$, 
  whence $r(\b V) \leq d-1$ follows.
  Furthermore, it is not at all obvious whether
  the ranks of the Moore-Penrose inverse $r( (\b C \wh {\b V}_N \b C' )^+ )$ 
  converge in probability to the rank $ r( (\b C \b V \b C')^+)$
  for a compatible contrast matrix $\b C$.
  Hence, the Wald-type statistic 
  $ N \wh {\b p}' \b C' (\b C \wh {\b V}_N \b C')^+ \b C \wh {\b p} $
  is not suitable for testing $H_0^p(\b C): \b C \b p = \b 0$:
  Its asymptotic behaviour is unclear and, hence, there is no reasonable choice of critical values.

  Instead, we utilize a statistic that does not rely on the uncertain convergence of ranks of generalized inverses.
  This leads us to the survival version of the so-called ANOVA-rank-type statistic
  \begin{equation}\label{eq:ATS}
  F_N(\b T) = \frac{N}{tr(\b T \wh {\b V}_N)} \wh {\b p}' \b T \wh {\b p},   
  \end{equation}
  where $\b T = \b C' (\b C \b C')^+ \b C$ is the unique projection matrix onto the column space of $\b C$. 
  Below we analyze both, its asymptotic behaviour under null hypotheses of the from $H_0^p(\b C): \b C \b p = \b 0$ and under the corresponding alternative hypotheses $H_a^p(\b C): \b C \b p \neq \b 0$.
  \begin{thm}
   \label{thm:teststat}
   Assume the asymptotic regime \eqref{eq:sample_size_conv} and that $tr(\b T {\b V})>0$.
   \begin{itemize}
    \item[a)] Under  $H_0^p(\b C)$ and as $N \rightarrow \infty$, we have
   $F_N(\b T) \oDo \chi = \b W' \b E_d' \b T \b E_d \b W/tr(\b T {\b V})$
   which is non-degenerate and non-negative with $E(\chi)=1$.
    \item[b)] Under $H_a^p(\b C)$ and as $N \rightarrow \infty$, we have
    $F_N(\b T) \oPo \infty.$
   \end{itemize}
  \end{thm}
  
  %
  %
  %
  As the distribution of $\chi$ depends on unknown quantities (cf. Theorem~\ref{thm:p}) the test statistic $F_N(\b T)$ in \eqref{eq:ATS} is no asymptotic pivot. To nevertheless 
  obtain proper critical values which lead to asymptotically exact inference procedures we next propose and study a resampling approach.

  \section{Inference via Multiplier Bootstrap}
  \label{sec:wbs}
  
  In this section, we apply suitably tailored multiplier bootstrap techniques
  in order to approximate the small sample distribution of $F_N( \b T )$.
  To this end, we consider the situation under $H_0^p(\b C)$ in which case we may expand
  \begin{align*}
   F_N( \b T ) & = \frac{N}{tr(\b T \wh {\b V}_N)} ( \wh {\b p} - \b p )' \b T ( \wh {\b p} - \b p ) 
     = \frac{N}{tr(\b T \wh {\b V}_N)} ( \d \phi_{\b S} \cdot ( \wh {\b S} -  \b S  ))' \b E_d' \ \b T \ \b E_d ( \d \phi_{\b S} \cdot ( \wh {\b S} - \b S )) + o_p(1) ,
  \end{align*}
  where $\wh {\b S}$ is the vectorial aggregation of all Kaplan-Meier estimators $\wh S_1, \dots, \wh S_d$.
  First, we replace the martingale residuals, that are attached to the Kaplan-Meier estimators, 
  with independent centered random variables which have approximately the same variance.
  In particular, we replace $\sqrt{N} (\wh {S}_i - S_i)$ with 
  $$\wh S(t) \cdot \sqrt{N} \sum_{k=1}^{n_i} G_{ik} \int_0^t [ (Y_i(u) - \Delta N_i (u)) Y_i(u) ]^{-1/2} \d N_{ik}(u); $$
  cf. \cite{dobler17} for a similar wild bootstrap Greenwood-type correction for tied survival data.
  Here we utilized the usual counting process notation \citep{abgk93}: 
  $N_{ik}$ indicates whether the event of interest already took place for individual $k$ in group $i$.
  The \emph{wild bootstrap} multipliers $G_{ik}, i=1, \dots, n_i, \ i=1, \dots, d$, 
  are i.i.d. with zero mean and unit variance and also independent of the data.
  In \cite{bluhmki18a, bluhmki18b} a similar multiplier resampling approach is
  applied to Nelson-Aalen and Aalen-Johansen estimators in one- and two-sample problems.
  
 In a next step toward the construction of a wild bootstrap statistic, we replace $\d \phi_{\b S}$ with $\d \phi_{\wh{\b S}}$.
 Let us denote the thus obtained wild bootstrap version of 
  $\sqrt{N} \d \phi_{\b S} \cdot ( \wh {\b S} - \b S )$ by
  ${\b W}_N^*$. Conditionally on the data, this $d^2$-variate random vector is for large $N$ 
  approximately normally distributed and its limit distribution coincides with that of $\b W$; see the proof of Theorem~\ref{thm:wbs} below for details.

 Finally, a wild bootstrap version $F_N^*(\b T)$ of $F_N(\b T)$ requires 
 that we also use a consistent wild bootstrap-type estimator $tr(\b T {\b V}_N^*)$ of $tr(\b T \wh {\b V}_N)$.
  It is found by replacing the estimators $\wh \Gamma_i$ with
  $$ \Gamma^*_i(s,t) = \wh S_i(s) \wh S_i(t) \ n_i \sum_{k=1}^{n_i} G_{ik}^2 \int_0^{s \wedge t} \frac{\d N_{ik}}{ (Y_i - \Delta N_i) Y_i}.$$
  Its conditional consistency was argued in \cite{dobler17} and a sufficient condition for this is $E(G_{11}^4) < \infty$.
  These wild bootstrap-type variance estimators also have the nice interpretation of optional variation processes of the wild bootstrapped Kaplan-Meier estimators (Dobler, 2017). 
  Hence, the resulting wild bootstrap version of $F_N( \b C )$ is
  \begin{align*}
  F_N^*( \b T ) & = \frac{1}{tr(\b T {\b V}_N^*)} {\b W}_N^{*}{\!}' \b E_d' \ \b T \ \b E_d {\b W}_N^*.
  \end{align*}
  The following conditional central limit theorem ensures the consistency of this resampling approach. 
  \begin{thm}
   \label{thm:wbs}
   Assume $E(G_{11}^4) < \infty$ and that the conditions of Theorem~\ref{thm:teststat} hold. 
   Conditionally on $(X_{ik}, \delta_{ik}), i =1, \dots, d, \ k = 1, \dots, n_i$,  we have for all underlying values of ${\b p}$
   $$F_N(\b T)^* \oDo \chi = {\b W}' \b E_d' \b T \b E_d {\b W}/ tr(\b T {\b V})$$
   in probability as $N \rightarrow \infty$.
  \end{thm}
 We would like to stress that the limit distribution coincides with that of $F_N(\b T)$ under $H_0^p(\b C)$. 
 For the wild bootstrap version $F_N(\b T)^*$, however, the convergence result holds under both, the null and the alternative hypothesis, 
 i.e. its conditional distribution always approximates the correct null distribution of the test statistic. 

  We conclude the theoretical part of this article with a presentation of deduced inference procedures for the effect sizes $\b p$. 
  To this end, let $c^*_{N,\alpha}$ denote the $(1-\alpha)$-quantile, $\alpha \in (0,1)$, of the conditional distribution of $F_N^*(\b T)$ given the data.
  In practice, this quantile is approximated via simulation by repeatedly generating sets of the wild bootstrap multipliers $G_{ki}$.
  \begin{cor}
   \label{cor:test}
   Under the assumptions of Theorem~\ref{thm:wbs}, the test
   $$ \varphi_N = 1\{ F_N(\b T) > c_{N,\alpha}^* \} $$
   is asymptotically exact and consistent.
   That is, $E( \varphi_N ) \rightarrow \alpha \cdot 1_{H_0^p(\b C)} + 1_{H_a^p(\b C)}$ as $N \rightarrow \infty$.
  \end{cor}
Now, let $r$ be the number of columns of $\b C'$ and denote by $\b c_1, \dots, \b c_r$ its column vectors.
The presentation of a simultaneous confidence region for the contrasts $\b c_\ell' \b p, \ell = 1, \dots, r,$ in Corollary~\ref{cor:sCIs} below will be done in an implicit manner. 
 \begin{cor}
   \label{cor:sCIs}
   Under the assumptions of Theorem~\ref{thm:wbs}, an asymptotically exact $(1-\alpha)$-confidence ellipsoid for the contrasts $\b c_\ell' \b p, \ell = 1, \dots, r,$   
   is given by 
   $$ CE = CE_{N,1-\alpha}(\b C) = \Big\{ \b v \in \R^r :
    (\b C \wh {\b p} - \b v)' (\b C \b C')^+ ( \b C \wh {\b p} - \b v) \leq \frac{tr(\b T \wh {\b V}_N)}{N} c_{N,\alpha}^* \Big\}.$$
   That is, $P( \b C \b p \in CE ) \rightarrow 1 - \alpha $ as $N \rightarrow \infty$.
  \end{cor} 
  
  
  \section{Simulations}\label{sec:simus}

  In this section, we assess the small sample properties of the test $\varphi_N$ as proposed in Corollary~\ref{cor:test}. 
  
  \subsection{Behaviour under null hypotheses}
  
\begin{figure}[ht]
\centering
 \includegraphics[width=0.8\textwidth]{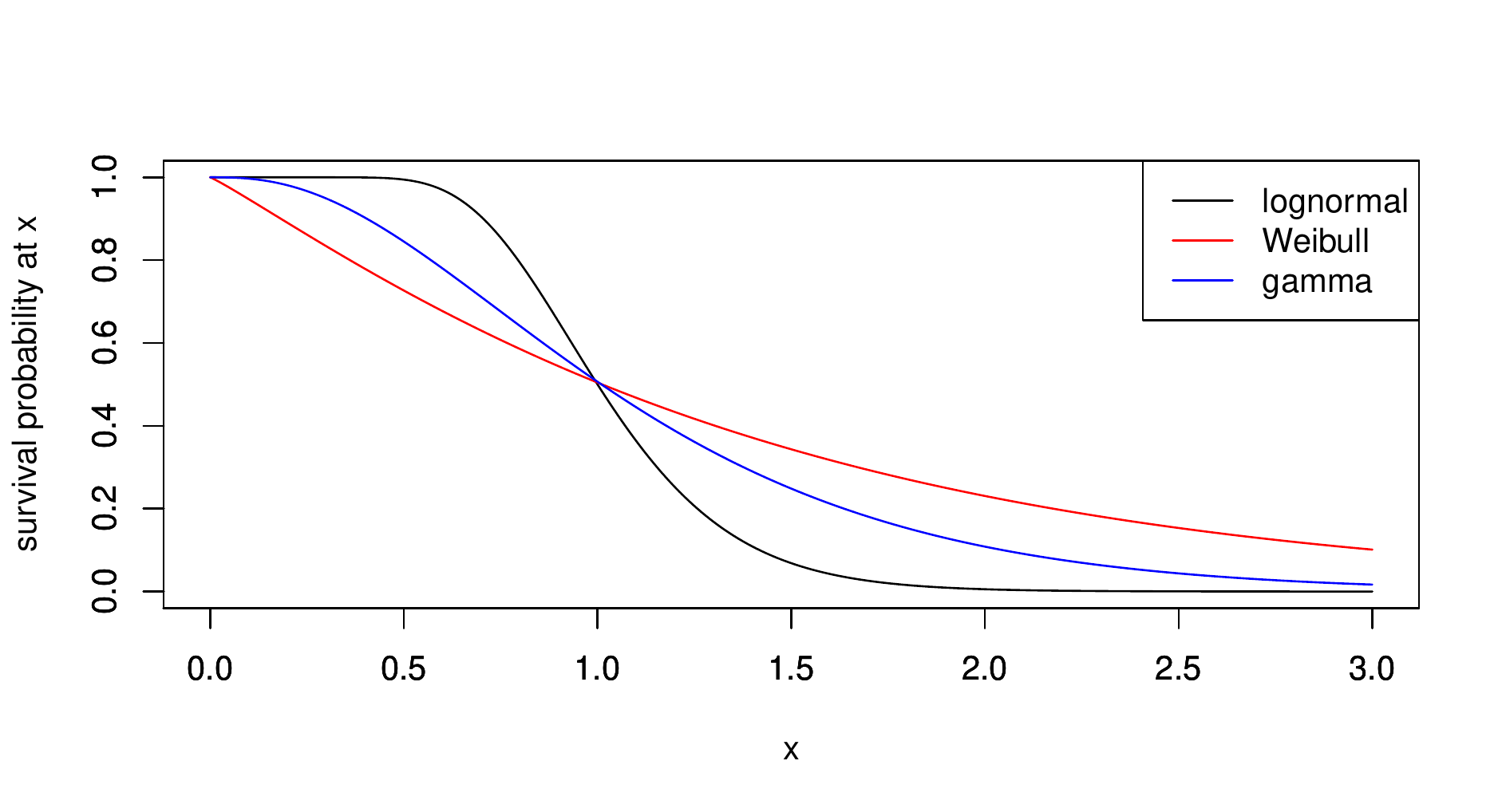}
 \caption{Survival functions underlying the first three simulated sample groups.}
 \label{fig:surv_func}
\end{figure}

    We first focus on its type I error control with respect to 
  \begin{itemize}
   \item various kinds of contrast matrices
   \item and different censoring intensities.
  \end{itemize}
{\bf Design and Sample Sizes.} For ease of presentation we restrict ourselves to a design with $d=6$ groups with different sample size layouts: 
  we considered small samples in a balanced design with $\b n_1 = (n_1,\dots,n_6)' = (10,10,10,10,10,10)'$ and two unbalanced designs with $\b n_2 = (n_1,\dots,n_6)' = (10,12,14,10,12,14)'$ and $\b n_3 = (10,12,14,14,10,12)'$, respectively. 
  To obtain designs with moderate to large sample sizes we increase these vectors component-wise by the factors $K \in \{2,3,5,10\}$. Moreover, depending on the question of interest, we below distinguish between a one-way layout with six independent groups and a $2\times 3$ two-way design.\\[0.4cm]
  {\bf Censoring Framework.} We considered exponentially distributed censoring random variables $C_{i1}\stackrel{ind}{\sim}Exp(\lambda_i)$ with the following vectors 
  $\bs \lambda = (\lambda_1,\dots,\lambda_6)'$ of rate parameters: 
    $\bs \lambda_1 = 0.4 \cdot \bs 1$, \ $\bs \lambda_2 = 0.5 \cdot \bs 1$, \ $\bs \lambda_3 = 2/3 \cdot \bs 1$, \
  $\bs \lambda_4 = (0.4, 0.5, 2/3, 0.4, 0.5, 2/3)'$, \ $\bs \lambda_5 = (0.4, 0.5, 2/3, 2/3, 0.5, 0.4)'$,
  where $\bs 1 \in \R^6$ is the vector consisting of $1$s only. 
  Thus, the first three settings correspond to equal censoring mechanisms with increased censoring rate from $\bs \lambda_1$ to $\bs \lambda_3$. The other two ($\bs \lambda_4$ and 
  $\bs \lambda_5$) lead to unequal censoring. By considering all $75$ possible combinations, many possible effects of censoring and sample size assignments are analyzed. 
  For example, in the set-up with $\b n_2$, $K=10$ and $\bs \lambda_4$, larger sample sizes are matched with a stronger censoring rate in an unbalanced design.\\[0.4cm]
  {\bf Contrast Matrices and Null Hypotheses.}
  We simulated the true significance level of the tests for the null hypotheses $H_0^p(\b C): \b C \b p = \b 0$ for two designs and different contrast matrices of interest:\\
  In case of a one-way design with $d=6$ groups we were interested in the null hypotheses of `{\it no group effect}'  or        
  `{\it equality of all treatment effects}' $H_0^p(\b C_1): \{ \b C_1 \b p = \b 0\} = \{p_1=\dots=p_6\}$. This may be described by considering the matrix $\b C_1 = \b P_6$, where here and below $\b P_d = \b I_d - \b J_d/d \equiv \b I_d - \bs 1_d \bs 1_d' /d$ denotes the $d$-dimensional centering matrix.\\
  Next, we consider a $2\times 3$ two-way layout with two factors $A$ (with two levels) and $B$ 
  (with three levels). This is incorporated in Model \eqref{eq:mod} by setting $d=a\cdot b= 2\cdot 3 =6$ and splitting up the index $i$ into two indices $i_1=1,2$ (for the levels of factor $A$) and $i_2=1,2,3$ (for the levels of factor $B$). Thus, we obtain survival times $T_{i_1i_2k}$, $k=1,\dots,n_{i_1i_2}$, and corresponding 
 nonparametric concordance effects $p_{i_1i_2}$. More complex factorial designs can be incorporated similarly. In this $2\times 3$ set-up we are now interested in testing the null hypotheses of 
  \begin{itemize}
   \item[(A)] `{\it No main effect of factor $A$}': $H_0^p(\b C_{2,A}): \{ \b C_{2,A} \b p = \b 0\} = \{ \bar{p}_{1\cdot} = \bar{p}_{2\cdot}\}$,
   \item[(B)] `{\it No main effect of factor $B$}': $H_0^p(\b C_{2,B}): \{ \b C_{2,B} \b p = \b 0\} = \{ \bar{p}_{\cdot 1} = \bar{p}_{\cdot 2} = \bar{p}_{\cdot 3}\}$ and
   \item[(AB)] `{\it No $A\times B$ interaction effect}': $H_0^p(\b C_{2,AB}): \{ \b C_{2,AB} \b p = \b 0 \} 
   = \{p_{i_1i_2} - \bar{p}_{i_1\cdot} -\bar{p}_{\cdot i_2} + \bar{p}_{\cdot\cdot} = 0 \text{ for all } i_1,i_2\}$,
  \end{itemize}
  where $\bar{p}_{i_1\cdot}$, $\bar{p}_{\cdot i_2}$ and $\bar{p}_{\cdot\cdot}$ denote the means over the dotted indices. In particular, the corresponding contrast matrices are given by 
   $ \b C_{2,A} = \b P_2 \otimes \frac13 \b J_3,
 \b C_{2,B} = \frac12 \b J_2  \otimes \b P_3,$ and $\b C_{2,AB} = \b P_2 \otimes \b P_3$,
where $\otimes$ indicates the Kronecker product.\\[0.4cm]
{\bf Survival Distributions.} For ease of presentation we only considered a rather challenging scenario, where the groups follow different survival distributions. In particular, we simulated 
\begin{itemize}
 \item[(G1)] a lognormal distribution with meanlog parameter $0$ and sdlog parameter $0.2726$ for the 
 first group,
 \item[(G2)]  a Weibull distribution with scale parameter $1.412$ and shape parameter $1.1$ for the second group,
 \item[(G3)] a Gamma-distribution with scale parameter $0.4$ and shape parameter $2.851$ for the third group and
 \item[(G4-G6)\!\!\!\!\!\!\!\!\!] \ \ \ \ \ \  mixing distributions of all pair combinations of the first three survival functions for the last three groups.
\end{itemize}
The first three survival functions are illustrated in Figure~\ref{fig:surv_func}. We note that preliminary simulations for more crude scenarios with identical survival distributions in all groups exhibited a much better type-$I$-error control of our testing procedure (results not shown).
Anyhow, the parameters of the above distributions were chosen in such a way that the nonparametric concordance effects of all groups are equal, i.e. $p_i = 0.5$ for all $i=1,\dots, 6$ (one-way) and 
$p_{i_1i_2} = 0.5$ for all $i_1=1,2; \ i_2=1,2,3$ (two-way), respectively. Thus, all considered null hypotheses are true. 
We would like to stress that the case of continuously distributed survival times corresponds to an infinite-dimensional problem and is thus more difficult than the discrete case. For example, the simulation study in \cite{dobler17} confirms this observation in a related problem: the convergence rate of the actual coverage probabilities of confidence bands 
to the nominal confidence level is much faster the more discretely distributed the survival data is. Moreover, to make the simulation scenario even more challenging, we considered the situation with infinite $\tau$ to also get an indication of the functionality of the test in this case.\\[0.4cm]
{\bf Simulations.}
  We chose as wild bootstrap multipliers centered unit Poisson variables because a formal Edgeworth expansion and two simulation studies in \cite{dobler17weird} indicate
  that  those have theoretical and practical advantages over the common choice of standard normal multipliers.
We chose the nominal level $\alpha=5\%$ and conducted each test 10,000 times for $K=1,2,3$ and 5,000 times for $K=5,10$ because of the massively increasing computational complexity for large samples.
Each test was based on critical values that were found using 1,999 wild bootstrap iterations. 
All simulations were conducted with the help of the R computing environment \citep{Rteam}.\\[0.4cm]
{\bf Results.}
The true type-$I$-error results for the four different null hypotheses are shown in Table~\ref{tab:level1} (left panel: one-way for $H_0^p(\b C_1)$ and right panel: two-way for 
$H_0^p(\b C_{2,A})$) and Table~\ref{tab:level2} (two-way for $H_0^p(\b C_{2,B})$ in the left and $H_0^p(\b C_{2,AB})$ in the right panel). 
It is apparent that all simulated  levels are elevated for the smallest sample sizes ($K=1$),
especially for the one-way test: here almost all type I error probabilities are between $13.0\%$ and $17.7\%$.
For the two-way tests, these probabilities are mainly between $8.1\%$ and $11.7\%$ in this case ($K=1$).
On the one hand, this is due to the relatively strong censoring rates:
for $\lambda = 0.4$, the censoring probabilities across all sample groups are between $33\%$ and $37\%$ (found by simulating 100,000 censoring and survival time random variables each), for $\lambda = 0.5$, these probabilities range from $39.5\%$ to $41.5\%$, and for $\lambda = 2/3$, they even reach values between $48.5\%$ and $49\%$; resulting in only $5$ to $7$ uncensored observations per group. On the other hand, not to restrict the time horizon in inferential procedures about survival functions appears to slightly slow down the convergence of type I error probabilities to the nominal level as the sample size increases; 
see \cite{dobler18} for similar findings in the context of confidence bands for unrestricted survival functions.
However, the error probabilities recover for samples of double size (i.e. between 20 and 28) already:
in the one-way design, these error rates drop to mainly $8.2\% - 9.9\%$, 
and in all two-way tests, we even achieve rates of mainly $6.1\% - 8\%$.
If the sample sizes are tripled (i.e. between 30 and 42), most of the type I error probabilities are between $7\%-8\%$ (one-way) or $5.2\%-6.9\%$ (two-way).
In case of the sample size factor $K=5$, all results are only slightly liberal,
and for $K=10$ (i.e. sample sizes between 100 and 140), we see that the nominal level is well attained.
              
\begin{table}[ht]
\centering
\begin{tabular}{cc|ccccc}
$n$ & \ \  $\lambda$  /  \ \ $K$ & 1 & 2 & 3 & 5 & 10 \\ \hline
            & $\bs \lambda_1$ & 14.7 & 8.8 & 7.2 & 6.4 & 5.7 \\ 
            & $\bs \lambda_2$ & 16.6 & 9.3 & 7.7 & 6.3 & 5.7 \\ 
  $\b n_1$  & $\bs \lambda_3$ & 19.7 & 11.0 & 8.6 & 6.6 & 5.8 \\ 
            & $\bs \lambda_4$ & 17.7 & 9.9 & 7.7 & 6.1 & 5.9 \\ 
            & $\bs \lambda_5$ & 17.4 & 9.5 & 7.8 & 6.7 & 6.3 \\ \hline
            & $\bs \lambda_1$ & 13.0 & 7.9 & 6.7 & 5.9 & 5.7 \\ 
            & $\bs \lambda_2$ & 13.9 & 8.4 & 7.1 & 6.9 & 5.4 \\ 
  $\b n_2$  & $\bs \lambda_3$ & 16.5 & 9.1 & 7.8 & 6.2 & 5.8 \\ 
            & $\bs \lambda_4$ & 14.9 & 8.5 & 7.3 & 6.0 & 5.6 \\ 
            & $\bs \lambda_5$ & 14.6 & 9.0 & 6.8 & 6.3 & 5.6 \\ \hline
            & $\bs \lambda_1$ & 12.2 & 8.2 & 7.1 & 6.0 & 5.2 \\ 
            & $\bs \lambda_2$ & 13.7 & 8.6 & 7.0 & 6.1 & 5.3 \\ 
  $\b n_3$  & $\bs \lambda_3$ & 17.7 & 9.3 & 7.6 & 6.8 & 5.9 \\ 
            & $\bs \lambda_4$ & 14.8 & 8.8 & 7.7 & 6.3 & 5.9 \\ 
            & $\bs \lambda_5$ & 14.1 & 8.6 & 7.2 & 6.3 & 5.9 \\ 
   \hline
\end{tabular}
\qquad
\begin{tabular}{cc|ccccc}
$n$ & \ \  $\lambda$  /  \ \ $K$ & 1 & 2 & 3 & 5 & 10 \\ \hline
            & $\bs \lambda_1$ & 9.0 & 6.3 & 6.0 & 5.7 & 5.4 \\ 
            & $\bs \lambda_2$ & 9.0 & 6.7 & 6.0 & 5.6 & 4.6 \\ 
  $\b n_1$  & $\bs \lambda_3$ & 10.5 & 7.0 & 6.3 & 5.9 & 5.6 \\ 
            & $\bs \lambda_4$ & 10.1 & 6.7 & 6.1 & 5.7 & 5.7 \\ 
            & $\bs \lambda_5$ & 9.4 & 6.2 & 5.8 & 5.6 & 5.0 \\ \hline
            & $\bs \lambda_1$ & 7.7 & 6.3 & 5.8 & 5.3 & 5.9 \\ 
            & $\bs \lambda_2$ & 8.1 & 6.3 & 6.0 & 6.0 & 5.3 \\ 
  $\b n_2$  & $\bs \lambda_3$ & 9.7 & 6.8 & 6.4 & 5.0 & 5.4 \\ 
            & $\bs \lambda_4$ & 8.2 & 6.3 & 6.2 & 5.9 & 5.4 \\ 
            & $\bs \lambda_5$ & 8.9 & 6.7 & 6.2 & 5.4 & 5.1 \\ \hline 
            & $\bs \lambda_1$ & 7.8 & 6.3 & 6.0 & 5.6 & 4.9 \\ 
            & $\bs \lambda_2$ & 8.5 & 6.2 & 5.5 & 5.3 & 5.0 \\ 
  $\b n_3$  & $\bs \lambda_3$ & 9.2 & 6.9 & 6.5 & 5.6 & 5.1 \\ 
            & $\bs \lambda_4$ & 8.4 & 6.1 & 6.1 & 5.8 & 4.5 \\ 
            & $\bs \lambda_5$ & 8.2 & 6.7 & 5.7 & 6.0 & 5.7 \\ 
   \hline
\end{tabular}
\caption{Simulated type I error probabilities in a one-way layout (left) and in a two-way design for main effect A (right) with sample size factor $K$.} 
\label{tab:level1}
\end{table}

\begin{table}[ht]
\centering
\begin{tabular}{cc|ccccc}
$n$ & \ \  $\lambda$  /  \ \ $K$ & 1 & 2 & 3 & 5 & 10 \\ \hline
            & $\bs \lambda_1$ & 10.0 & 7.2 & 6.4 & 6.2 & 6.0 \\ 
            & $\bs \lambda_2$ & 11.4 & 7.7 & 6.7 & 5.9 & 4.9 \\ 
  $\b n_1$  & $\bs \lambda_3$ & 13.4 & 8.0 & 6.9 & 5.9 & 5.8 \\ 
            & $\bs \lambda_4$ & 12.2 & 7.6 & 6.9 & 6.1 & 5.9 \\ 
            & $\bs \lambda_5$ & 12.1 & 7.5 & 6.7 & 5.9 & 5.6 \\ \hline
            & $\bs \lambda_1$ & 9.5 & 6.6 & 6.1 & 6.0 & 5.0 \\ 
            & $\bs \lambda_2$ & 10.2 & 7.4 & 6.5 & 6.0 & 5.5 \\ 
  $\b n_2$  & $\bs \lambda_3$ & 11.6 & 7.8 & 6.6 & 5.6 & 5.7 \\ 
            & $\bs \lambda_4$ & 10.4 & 7.0 & 6.4 & 5.5 & 6.1 \\ 
            & $\bs \lambda_5$ & 10.2 & 7.1 & 6.2 & 5.9 & 5.4 \\ \hline
            & $\bs \lambda_1$ & 9.5 & 7.2 & 5.8 & 5.2 & 5.2 \\ 
            & $\bs \lambda_2$ & 9.6 & 6.8 & 6.3 & 5.4 & 5.6 \\ 
  $\b n_3$  & $\bs \lambda_3$ & 11.6 & 7.4 & 6.9 & 6.2 & 5.6 \\ 
            & $\bs \lambda_4$ & 9.9 & 7.4 & 6.2 & 6.5 & 5.4 \\ 
            & $\bs \lambda_5$ & 9.7 & 7.2 & 6.4 & 5.7 & 5.0 \\ 
   \hline
\end{tabular}
\qquad
\begin{tabular}{cc|ccccc}
$n$ & \ \  $\lambda$  /  \ \ $K$ & 1 & 2 & 3 & 5 & 10 \\ \hline
            & $\bs \lambda_1$ & 10.1 & 7.2 & 6.3 & 5.7 & 5.3 \\ 
            & $\bs \lambda_2$ & 11.2 & 7.2 & 6.2 & 5.9 & 5.1 \\ 
  $\b n_1$  & $\bs \lambda_3$ & 13.3 & 8.5 & 7.0 & 6.5 & 5.5 \\ 
            & $\bs \lambda_4$ & 11.6 & 7.8 & 6.6 & 6.1 & 5.3 \\ 
            & $\bs \lambda_5$ & 11.6 & 7.7 & 6.4 & 5.9 & 5.6 \\  \hline
            & $\bs \lambda_1$ & 9.2 & 6.6 & 5.9 & 6.2 & 5.3 \\ 
            & $\bs \lambda_2$ & 9.8 & 6.9 & 6.6 & 5.4 & 5.7 \\ 
  $\b n_2$  & $\bs \lambda_3$ & 11.7 & 7.5 & 6.4 & 5.8 & 5.4 \\ 
            & $\bs \lambda_4$ & 9.8 & 6.8 & 6.4 & 5.2 & 5.0 \\ 
            & $\bs \lambda_5$ & 10.8 & 7.1 & 6.2 & 5.8 & 5.6 \\  \hline
            & $\bs \lambda_1$ & 8.3 & 6.6 & 6.3 & 5.2 & 5.1 \\ 
            & $\bs \lambda_2$ & 9.6 & 6.9 & 5.9 & 5.8 & 5.7 \\ 
  $\b n_3$  & $\bs \lambda_3$ & 11.2 & 7.6 & 6.4 & 5.3 & 5.8 \\ 
            & $\bs \lambda_4$ & 10.4 & 6.9 & 5.8 & 5.4 & 4.7 \\ 
            & $\bs \lambda_5$ & 9.8 & 6.8 & 5.5 & 5.5 & 5.7 \\ 
   \hline
\end{tabular}
\caption{Simulated type I error probabilities in a two-way design for main effect B (left) and for interaction effect AB (right) with sample size factor $K$.} 
\label{tab:level2}
\end{table}

\subsection{Behaviour under shift alternatives}

\begin{figure}[ht]
   \includegraphics[width=\textwidth]{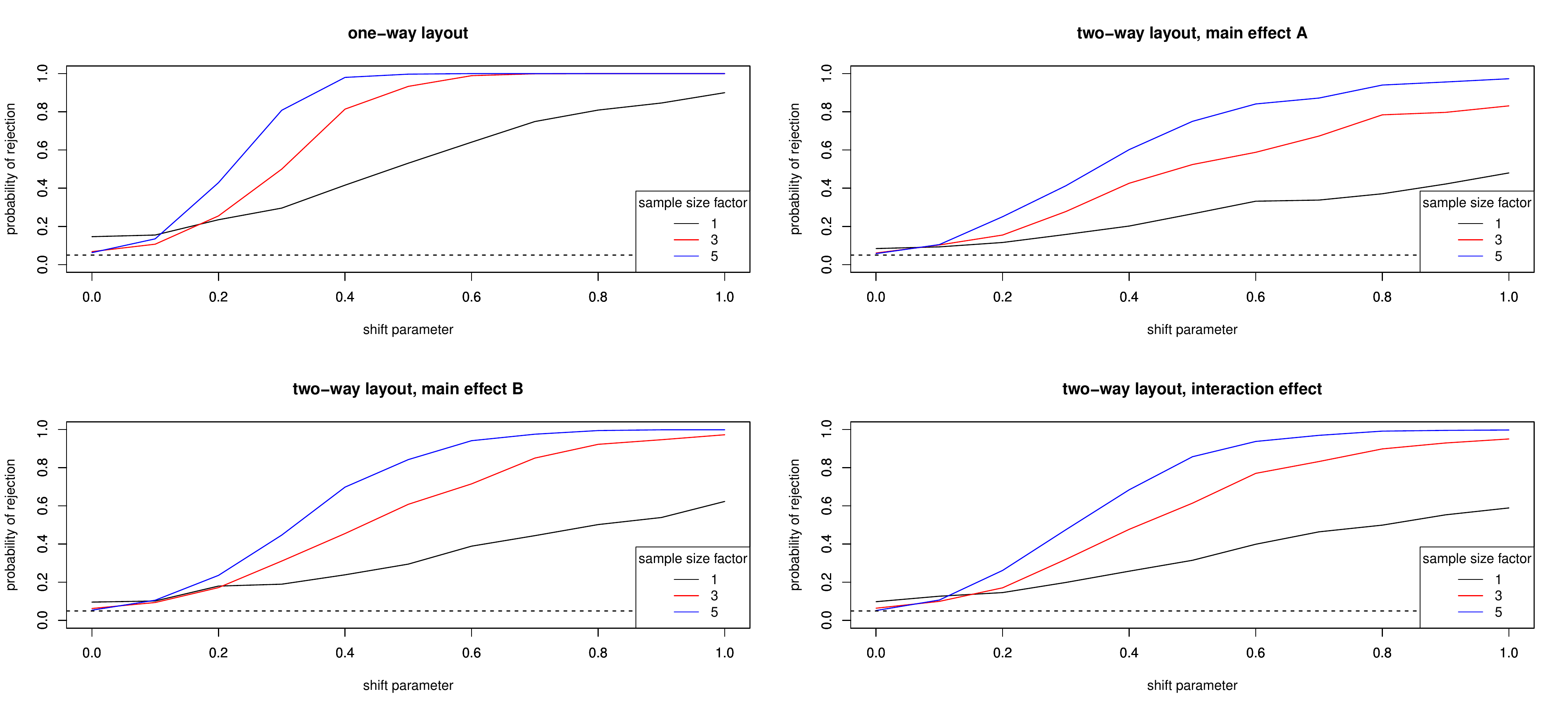}
   \caption{Power functions for shift alternatives for different null hypotheses: in the one-way layout  (sample sizes $\b{n} = K \cdot \b{n}_2$, censoring rates $\bs \lambda = \bs \lambda_5$), in the two-way layout for main effect A ($\b{n} = K \cdot \b{n}_3$, $\bs \lambda = \bs \lambda_4$), for main effect B ($\b{n} = K \cdot \b{n}_3$, $\bs \lambda = \bs \lambda_2$), and for the interaction effect ($\b{n} = K \cdot \b{n}_2$, $\bs \lambda= \bs \lambda_4$), $K = 1,3,5$. The nominal significance level is $\alpha =5\%$ (- - -).}
   \label{fig:power}
  \end{figure}

In addition to the simulations of the previous subsection, we also conducted a small power simulation of the above tests. For the alternative hypotheses, we considered a shift model: taking the same six basic survival and censoring functions as in the first set of simulations, we shift all survival and censoring times of the first sample group by $\delta \in \{0.1, 0.2, \dots, 1\}$.
In this way, we maintain the same censoring rates as before and the distance to the null hypotheses is gradually increased: for growing $\delta >0$, we obtain a growing relative effect $p_1 > 0.5$ (one-way) and 
$p_{11} > 0.5$ (two-way), respectively. For each of the above considered contract matrices, $\b C_1, \b C_{2,A}, \b C_{2,B}, \b C_{2,AB}$,
we conducted one set of simulations with different unbalanced sample sizes and censoring rate combinations. 
For each set-up, we increased the sample sizes by the factors $K=1,3,5$. The results are displayed in Figure~\ref{fig:power}.

We see that, even for the smallest sample sizes (between 10 and 14), the power of the two-way testing procedures increase to $0.5$ or $0.6$ as the shift parameter approaches 1. For larger samples sizes the theoretically proven consistency is apparent. In comparison, the one-way test has a much higher power: 
For the undersized case ($K=1$) it already reaches a power of $0.8$ while for moderate to larger sample sizes the power is almost $1$ for shift parameters $\delta\geq 0.5$. In comparison to the two-way procedure its superior power is, however, partially paid at the price of its pronounced liberality; especially for small sample sizes.

All in all, the simulations confirm that all tests have a satisfactory power with increasing sample size and/or shift parameter while maintaining a reasonable control of the nominal level for sample sizes of $30$ to $42$ already.

  \section{Data example}\label{sec:data_Example}
  
\begin{figure}
   \includegraphics[width=\textwidth]{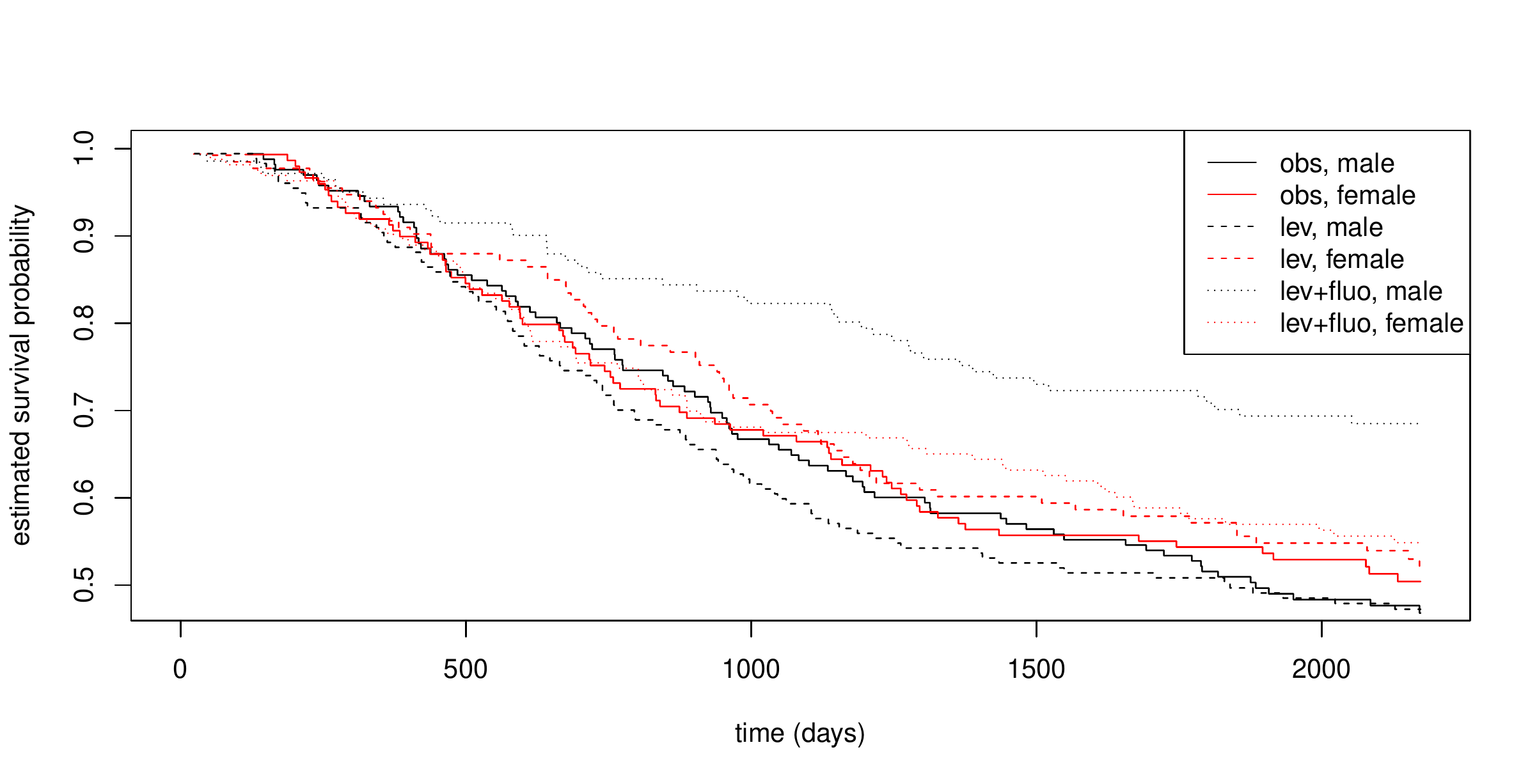}
   \caption{Kaplan-Meier estimator for male and female subgroups, discriminated further according to treatment: obs = observation, lev = levamisole treament, lev+fluo = combined levamisole and fluorouracil treatment. The end time in the plot is $\tau$ = day  2173.}
   \label{fig:kmes}
  \end{figure}

  We illustrate the developed theory on a dataset from a colon cancer study \citep{moertel90}. 
  Considering the patients in \emph{Stage C}, that is, there had been metastases to regional lymph nodes,
  the data consist of eligible 929 patients suffering from colon cancer. Survival (measured in days) was the primary endpoint of the study. We focus on the two factors `gender' and `treatment' (with three levels) to obtain a crossed $2\times 3$ survival design which is in line with a setting from our simulation study. In particular, there were
  315 patients in the observation group,
  310 others were treated with levamisole,
  and 304 received levamisole, combined with fluorouracil.
  Levamisole was originally used as an anthelmintic drug and fluorouracil (5-FU) is a medicine to treat various types of cancer. 
  The patients in the study had been randomized into one of these three treatment groups.
  Also, there were nearly as many women (445) as men (484) involved in the study.
  Figure~\ref{fig:kmes} depicts the Kaplan-Meier estimates of the survival probabilities for each treatment $\times$ sex subgroup.
  We refer to \cite{moertel90} for more details about the study.
  The dataset is freely accessible via the \texttt{R} command \texttt{data(colonCS)} after having loaded the package \emph{condSURV} \citep{meira16,meira16R}. 
  
  The aim is now to investigate the presence of main or interaction effects of treatment and gender. As there are several ties in the data (roughly $16\%$; see Appendix~\ref{app:add_info} for details) and we do not want to impose specific distributional assumptions, we focus on the nonparametric concordance effects. To this end, we first have to choose a proper $\tau$. From our retrospective view, the most reasonable choice is found by determining for each group the minimal observed censoring time that exceeds all observed survival times in that group. We call these censoring times ``terminal times''. Then,  $\tau$ is set to be the minimal terminal time.   In doing so, the group with that minimum terminal time does not benefit nor does it suffer from having the earliest terminal time when compared to the other groups.
  
  The first block in Table~\ref{tab:dataex_descr} shows the sample sizes of all subgroups.
  In the present data example, the minimal terminal time is $\tau = 2173$; see the second block of Table~\ref{tab:dataex_descr}.
  In view of the sample sizes and the censoring rates given in the third block of Table~\ref{tab:dataex_descr},
  we compare the present dataset with the simulation set-ups in Section~\ref{sec:simus}: a
  similarly strong censorship is obtained for $\bs\lambda_3$ and comparable sample sizes $n \in [100,140]$ for the choice $K=10$. Thus, judging from the rightmost columns of Tables~\ref{tab:level1} and~\ref{tab:level2}, we find it is safe to assume actual type I error probabilities of about $5.1\%$ to $5.9\%$ of the proposed nonparametric one- and two-way survival tests. 

  \begin{table}[ht]
   \centering
   \begin{tabular}{l|cc|cc|cc|cc}
    & \multicolumn{2}{|c}{sample size} & \multicolumn{2}{|c}{terminal time} & \multicolumn{2}{|c}{censoring rate} & \multicolumn{2}{|c}{effect size}\\
    treatment & male & female & male & female & male & female & male & female \\ \hline 
    observation & 166 & 149 & 2800 & 2562 & 47.6 & 51.0 & 0.475 & 0.483  \\
    levamisole & 177 & 133 & 2915 & 2173 & 47.5 & 52.6 & 0.459 & 0.501 \\
    levamisole plus fluorouracil & 141 & 163 & 2726 & 2198 & 68.8 & 55.2 & 0.581 & 0.501 \\ \hline
   \end{tabular}
    \caption{For each subgroup: sample size, smallest censoring time (in days) exceeding the largest survival time, censoring rate (in \%) after taking the minimum of each event time and $\tau = 2173$, and nonparametric concordance effects. Columns: sex, row: treatments.}
    \label{tab:dataex_descr}
  \end{table}
  
  We tested the data in one- and two-factorial set-ups and chose $\alpha = 5\%$ as the significance level.
  As in the simulation study, we used $B=1,\!999$ bootstrap iterations for each test.
  For the tests in the two-factorial model, we considered the null hypotheses corresponding to no main treatment effect, no main effect in sex, and no interaction effect between both.
  The test results, by means of p-values, are shown in Table~\ref{tab:dataex_pval}.
  
  \begin{table}[ht]
   \centering
   \begin{tabular}{lll|c}
    Null hypothesis & $H_0^p ( \cdot ) $ & set-up &  p-value \\ \hline 
    Equality of all effects & $H_0^p (\b C_1)$ & one-factorial & $ < 0.001$ \\
    No main effect in sex & $H_0^p(\b C_{2,A})$ & two-factorial  & $ 0.331 $ \\
    No main effect in treament & $H_0^p(\b C_{2,B})$ & two-factorial & $< 0.001 $ \\
    No interaction effect & $H_0^p(\b C_{2,AB})$ & two-factorial & $< 0.001 $ \\ \hline
   \end{tabular}
    \caption{p-values of different hypothesis tests for the anaylsis of the colonCS data-set.}
    \label{tab:dataex_pval}
  \end{table}

  We found a significant indication against the equality of all $d=6$ groups (p-value $ < 0.001$) but this difference between groups could not be inferred to result from a difference between the sexes (p-value $= 0.331$).
  However, we found a significant treatment effect (p-value $ < 0.001$)
  as well as a significant interaction effect between treatment and sex (p-value $ < 0.001$).
  The $p$-values in Table~\ref{tab:dataex_pval} have not been adjusted for a type I error multiplicity but it is obvious that the results remain the same after an application of, say, the Bonferroni-Holm procedure.
  
  Indeed, looking at the rightmost block of Table~\ref{tab:dataex_descr},
  we agree with the findings of the hypothesis tests:
  the gender effect seems to be cancelled out if the treatment groups are combined,
  but within the male gender there seems to be a big difference in the concordance effects ($p_{1i_2} \in [0.459, 0.581]$).
  Also, the interaction effect is apparent,
  as the female groups do not seem to strongly benefit from any treatment ($p_{2i_2} \in [0.483, 0.501]$).
  On the other hand, the male groups exhibit a worse than average survival fitness in the observation and the levamisole treatment group ($p_{11}=0.475,\ p_{12}=0.459$)
  but a much better than average survival fitness for the combination treatment ($p_{13} = 0.581$). Here the value 
  $p_{13} = 0.581$ roughly means that a randomly generated observation from this specific group survives 
  a randomly generated observation from the mean distribution of all groups with probability $58.1\%$. Taking another look at the Kaplan-Meier curves in Figure~\ref{fig:kmes},
  we immediately see that our concordance effects and the test outcomes make sense.
  We clearly see that there is a big difference in the male survival probabilities
  (the combination treatment group is superior to the levamisole treatment group which is in turn superior to the observation group)
  but there is not much of a difference between the female groups' survival curves.
  Indeed, comparing the Kaplan-Meier curve of the pooled males' survival times with that of the pooled females' times, 
  we graphically find no evident main gender effect.
  The plot of both Kaplan-Meier estimators is shown in Appendix~\ref{app:add_info}.
  
  Finally, we relate our results to the original findings of \cite{moertel90}
  whose analyses involve the Cox proportional hazards model and logrank tests.
  The authors also detected that ``Therapy with levamisole plus fluorouracil produced an unequivocal advantage over observation'' and that levamisole alone did not produce a detectable effect.
  Furthermore, they concluded from an exploratory subset analysis that the ``levamisole-fluorouracil treatment appeared to have the greatest advantage among male patients [...]''.
  This is exactly what we confirm in our analysis based on the nonparametric two-factorial tests.
  However, \cite{moertel90} neither account for the present ties in the data (note that the Cox regression postulates continuous outcomes) nor did they clearly stress the rather weak effect of the  levamisole-fluorouracil treatment for women. 
  They just state that their ``results show [...] striking contradictions to those of subset analyses reported in the NCCTG study, in which levamisole plus fluorouracil was found to be most effective in reducing the risk of recurrence among female patients [...]'' among other subgroups of patients.

  \section{Discussion}\label{sec:dis}
  
We proposed novel nonparametric inference procedures for the analysis of factorial survival data that may be subject to independent random right-censoring. Critical values are obtained from a multiplier wild bootstrap approach leading to 
asymptotically valid tests and confidence regions for meaningful effect parameters. Thereby, the procedures do not require any multiplicative or additive hazard structure nor specific distributional survival and censoring assumptions.  In particular, different group distributions are allowed and 
ties are accounted for accordingly. Moreover, different to the nonparametric survival procedures of \cite{akritas97} and \cite{akritas2011nonparametric}, our methods are not only driven towards hypothesis testing but also to uncertainty  quantification of the underlying effect estimators. The latter can be used to comprehensibly describe and infer main and interaction effects in general nonparametric factorial survival designs with an arbitrary number of fixed factors. 
Together with a $1$-$1$ connection with hazard ratios in proportional two-sample designs \citep{bruckner2017sequential},
this makes the new methods appealing for practical purposes. 

To investigate their theoretical properties, we rigorously proved central limit theorems of the underlying statistics and consistency of the corresponding procedures. In addition, 
extensive simulations were conducted for one- and two-way designs to also assess their finite sample properties in terms of power and type-$I$-error control. In case of small sample sizes with less than $10$ completely observed subjects per group, they revealed a liberal behaviour; especially for the one-way testing procedure. However, for moderate to larger sample sizes the asymptotic results kicked in and the stated theoretical results were recovered.

Finally, the methods were used to exemplify the analysis of survival data in a study about treatments for colon cancer patient within a two-factorial survival design. As severe ties were present in the data, classical hazard based methods were not directly applicable. In comparison, our newly proposed nonparametric methods provided a very decent alternative for the analysis of such factorial survival designs without postulating any strict assumptions.

To allow for a straightforward application, it is planned to implement the procedure into an easy to use R-package. 
In future research we will consider the case of stochastically ordered subgroups, for which a multiple testing algorithm could be developed with the aim to detect significantly different collections of all subgroups:
subgroups with no significant differences in the nonparametric concordance effects may be combined to facilitate the interpretation of the outcomes
and to ultimately serve for the development of different, more personalized medicines, one for each new subgroup combination.
Moreover, extensions of the current methodology to ordered alternatives or factorial designs obtained via stratified sampling will be part of a practically useful consecutive testing procedure.

  \section*{Acknowledgements}
   Markus Pauly likes to thank for the support from the German Research Foundation (Deutsche Forschungsgemeinschaft).

   \bibliographystyle{plainnat}
  \bibliography{literatur}

  \appendix
  
  \section*{Appendix}
  
  \section{Proofs}
  \label{app:proofs}
  
  \subsection*{Proof of Lemma~\ref{lem:w}}
  
  The proof of this lemma is based on an extension of the results of \cite{dobler2016bootstrap}.
  In that article it is shown that the functional
  $\widetilde \phi : (BV[0,\tau])^2 \rightarrow \R, \ (f,g) \mapsto - \int_0^\tau f^\pm(u) \d g(u)$
  is Hadamard-differentiable in functions of bounded variation; see the proof of their Theorem~1 in the supplement to that article for more details.
  Hence, a similar result transfers to the functional $\phi$ in the present article because $\phi$ is simply a multivariate aggregation of functionals of the same type as $\widetilde \phi$.
  
  For the sake of completeness, we point out that the $((i-1)d+j)$th entry of continuous and linear Hadamard-derivative of $\phi$ at $\b S$ is given by
 $$(\d \phi_{\b S} \cdot \b h)_{i,j} = \int h_j^\pm \d S_i - \int h_i^\pm \d S_j, \quad \b h = (h_1, \dots, h_d)' \in (D[0,\tau])^d,$$
 if $h_1(0) = \dots = h_d(0) = 0$ which indeed is the case in the present application. 
  
  The other fundamental result that we require is the convergence of all normalized Kaplan-Meier estimators in distribution.
  It is well-known that $\sqrt{n_i} ( \wh S_i - S_i) \oDo U_i $ on the Skorokhod space $D[0,\tau]$ as $n_i \rightarrow \infty$ in the present situation of independent right-censoring.
  Hence, due to the independence of all sample groups, it also follows that
  $diag(n_1, \dots, n_d) \cdot (\wh {\b S} - \b S) \oDo \b U $ on $(D[0,\tau])^d$ as $\min(n_1, \dots, n_d) \rightarrow \infty$.
  
  Thus, after having added the multiplicative term $\frac{\sqrt{n_i}}{\sqrt{n_i}}$ in each component $i=1, \dots, d$,
  it follows by means of the functional delta-method (cf. Theorem~3.9.4 in \citealp{vaart96}) 
  that $\sqrt{N} ( \phi(\wh {\b S}) - \phi(\b S))$ converges in distribution as asserted.
  \qed
  
  \subsection*{Proof of Theorem~\ref{thm:p}}
  
  This convergence result follows immediately from Lemma~\ref{lem:w} in combination with the continuous mapping theorem.
  \qed
  
  \subsection*{Proof of Lemma~\ref{lem:gamma}}

  It is well-known that each Kaplan-Meier estimator $\wh S_i$
  and also each covariance function estimator $\wh \Gamma_i$ is uniformly consistent on $[0,\tau]$ for $S_i$ and on $[0,\tau]^2$ for $\Gamma_i$, respectively, $i=1, \dots, d$.
  Note here that due to the assumption $P(T_{i1} > \tau) > 0$ for all groups $i$ the cumulative hazard functions are bounded on the interval of interest: $\Lambda_i(\tau) < \infty$.
  
  Since $\wh {\b V}_N$ is a continuous functional of these estimators,
  it follows from the continuous mapping theorem that the matrix $\wh {\b V}_N$ is consistent for $\b V$.
  Likewise, we conclude that the estimator $tr(\b T \wh {\b V}_N)$ is consistent for the trace $tr(\b T {\b V})$ due to continuity of the trace and linear maps. \qed
  
   \subsection*{Proof of Theorem~\ref{thm:teststat}}
  
  Again, the result essentially follows from the continuous mapping theorem as $\b p \mapsto \b p' \b T \b p$ is a continuous map.
  Thus, taking Lemma~\ref{lem:gamma} into account, an application of Slutsky's lemma reveals that
  $$F_N(\b C) \stackrel{d}{\longrightarrow} \chi$$
  under $H_0^p(\b C)$ as $N \rightarrow \infty$.
  
  To show that $\chi$ has unit expectation,
  we note that the trace in the denominator of $F_N(\b C)$ accounts for the asymptotic expectation of the quadratic form term:
  the expectation of a quadratic form $\b X' \b A \b X$ for a random vector $\b X$ with mean $\bs \mu = E(\b X)$ and covariance matrix $\bs \Sigma = cov( \b X )$ equals
  $ E(\b X' \b A \b X ) = tr( \b A \bs \Sigma ) + \bs \mu' \b A \bs \mu$;
  cf. Corollary~3.2b.1 in  \cite{mathai92}.
  In our situation, where $\b A = \b T$ and $\b T \b p = \b 0$ under $H_0^p(\b C)$,
  the asymptotic expectation of $N \wh {\b p}' \b T \wh {\b p}$ under $H_0^p$ thus simplifies to $tr( \b T \bs V ) $.

  On the other hand, under the alternative hypothesis $H_a^p(\b C)$, we have that 
  $$ F_N(\b T) =  \frac{N}{tr(\b T \wh {\b V}_N)} \wh {\b p}' \b T \wh {\b p} =  \frac{N}{tr(\b T \wh {\b V}_N)} (\wh {\b p} - \b p)' \b T (\wh {\b p} - \b p) + \frac{N}{tr(\b T \wh {\b V}_N)} (2 \wh {\b p} - \b p)' \b T \b p. $$
  The first term on the right-hand side converges in distribution to $\chi$ as $N \rightarrow \infty$.
  Of the second term, $(2 \wh {\b p} - \b p)' \b T \b p$ converges to the positive number $\b p' \b T \b p$.
  Thus, an application of Slutzky's lemma yields that $F_N(\b T) \stackrel{p}{\rightarrow} \infty$ as $\min \bs n \rightarrow \infty$.
  \qed
  
   \subsection*{Proof of Theorem~\ref{thm:wbs}}
  
  This conditional central limit theorem holds
  because for each wild bootstrapped Kaplan-Meier estimator a conditional central limit theorem holds separately; cf. \cite{dobler17} who suggested a wild bootstrap resampling technique for independently right-censored competing risks and survival data in the presence of ties.
  They are combined with the help of the continuous mapping theorem 
  and also with the consistency of the wild bootstrapped covariance estimators which were shown in that same paper.
  Hence, the continuous mapping theorem implies that the conditional distribution of $ F_N^*(\b T)$ given the data converges weakly in probability to the distribution of $\chi$.
  \qed
  
  \section{Asymptotic Covariance Matrix $\b V$ of $\sqrt N ( \wh {\b p} - \b p)$}
  \label{app:asy.cov}
  
  The covariance matrix $\b V$ in Theorem~\ref{thm:p} has the compact form
  \begin{align*}
   \b V = &  \frac1{d^2} \sum_{j=1}^d \frac1{\lambda_j} \int \int \Gamma_j^{\pm \pm} \d \b S \d \b S'
    + \int \int diag \Big( \Big(\frac1{\lambda_i} \Gamma_i^{\pm \pm} \Big)_{i=1}^d \Big) \d \bar S \d \bar S \\
     & - \frac1d \int \int \Big(\frac1{\lambda_i} \Gamma_i^{\pm \pm} \Big)_{i=1}^d \d \bar S \d \b S'
    - \Big( \frac1d \int \int \Big(\frac1{\lambda_i} \Gamma_i^{\pm \pm} \Big)_{i=1}^d \d \bar S \d \b S' \Big)' \\
   = & \frac{1}{d^2}  \int \int \b 1_d' diag(\bs \lambda)^{-1} \bs \Gamma^{\pm \pm} \d \b S \d \b S'
    + diag(\bs \lambda)^{-1} \int \int diag(\bs \Gamma^{\pm \pm}) \d \bar S \d \bar S \\
    & - \frac1d \int \int diag(\bs \lambda)^{-1} \bs \Gamma^{\pm \pm} \d \bar S \d \b S'
    - \Big( \frac1d \int \int diag(\bs \lambda)^{-1} \bs \Gamma^{\pm \pm} \d \bar S \d \b S' \Big)' ,
  \end{align*}
  where $\bs \Gamma = (\Gamma_i)_{i = 1, \dots, d}$ 
  is the vectorial aggregation of asymptotic covariance functions of all Kaplan-Meier estimators.

  \section{Additional Information on the Data Analysis}
  \label{app:add_info}
  As stated in Section~6 of the paper there are severe ties present in the dataset. The explicit numbers of ties are presented in Table~\ref{tab:ties1} below. They show that $130$ different observations occurred at least twice in the data set. In particular, only $781$ different realizations are encountered among all $929$ individuals leading to 
a ties rate of roughly $148/929\approx 16\%$.              
\begin{table}[ht]
\centering
\begin{tabular}{|c|c|c|c|}
\hline\multicolumn{4}{|c|}{Observations that appear}\\
once & twice & three times & four times\\\hline
651 &115 & 12 & 3\\\hline
\end{tabular}
\caption{Number of ties in the data before we truncated at $\tau=2173$ days} 
\label{tab:ties1}
\end{table}

\noindent After choosing the truncation time $\tau=2173$ days as described in the paper, the numbers of ties further increase as illustrated in Table~\ref{tab:ties2}.
This leads to a ties rate of approximately $439/929\approx 47\%$ in the final data analysis.

\begin{table}[ht]
\centering
\begin{tabular}{|c|c|c|c|c|}
\hline\multicolumn{5}{|c|}{Observations that appear}\\
once & twice & three times & four times& 361 times\\\hline
421 & 59 & 7 & 2 & 1\\\hline
\end{tabular}
\caption{Number of ties in the data after we truncated at $\tau=2173$ days} 
\label{tab:ties2}
\end{table}

\begin{figure}[hb]
 \centering
 \includegraphics[width=0.63\textwidth]{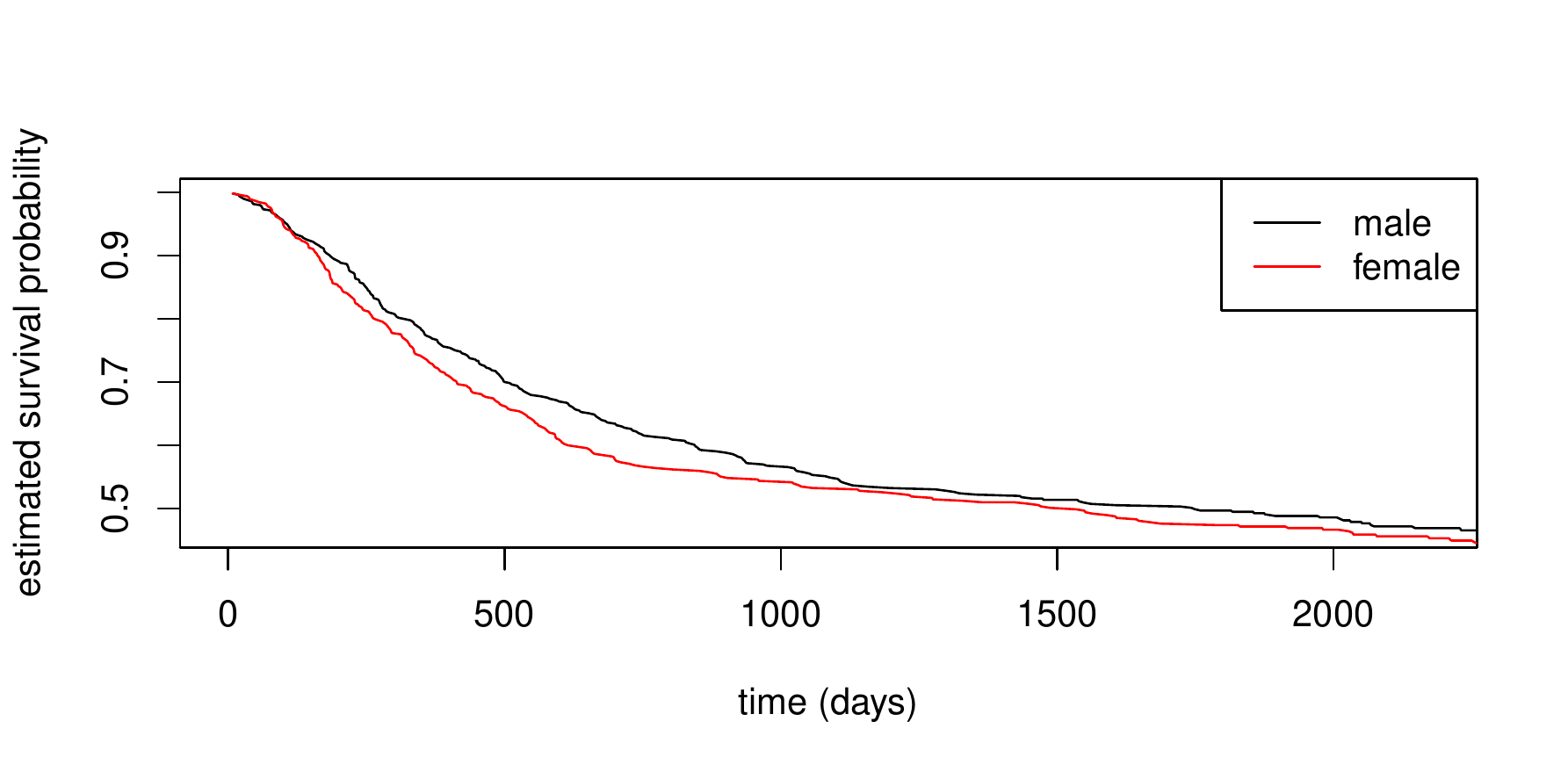}
 \caption{Kaplan-Meier curves for male (black) and female patients (red).}
 \label{fig:kmes_male_female}
\end{figure}

\noindent Finally, we present in Figure~\ref{fig:kmes_male_female} the Kaplan-Meier curves with respect to gender, thus after a combination of all three male subgroups and all three female subgroups.
Even though the test for the main effect in gender was not detected to be significant in the set-up with six subgroups, we found that gender became highly significant if the two-sample problem is considered (p-value $ < 0.001$).
Indeed, in Figure~\ref{fig:kmes_male_female} the survival times of male patients appear to be slightly stochastically larger than those of the females.
We thus see that the tests further gain in power when different subgroups are combined as the involved sample sizes increase a lot. 
Of course, one should be careful with this new result because combining subgroups may not always make sense. In this case, especially in view of the significant interaction effect between treatment and gender, it could be reasonable to treat women differently than men.

\end{document}